\documentclass[12pt,fleqn]{amsart}

\usepackage[all]{xy}
\usepackage{tikz-cd}
\usepackage{amscd}
\usepackage{setspace}
\usepackage{korpi-lib}

\usepackage{comment}
\usepackage{enumerate}
\usepackage{color}
\usepackage{enumerate}
\usepackage{amssymb}
\usepackage{hyperref}
\newcommand{\bmd}{d_{\rm BM}}

\numberwithin{equation}{section}

\title {How many miles from $L_\infty$ to $\ell_\infty$?}
\author[M. Korpalski]{Maciej Korpalski}
\author[G.\ Plebanek]{Grzegorz Plebanek}
\address{Instytut Matematyczny, Uniwersytet Wroc\l awski, pl.\ Grunwaldzki 2, 50-384 Wroc\-\l aw\\ Poland}
\email{Maciej.Korpalski@math.uni.wroc.pl \\ Grzegorz.Plebanek@math.uni.wroc.pl}
\date{}

\thanks{The first author was partially supported by {\em Młody Badacz}  grants funded by University of Wrocław}

\begin{document}

\begin{abstract}
The classical Banach spaces $L_\infty[0,1]$ and $\ell_\infty$ are isomorphic. We present here some
lower and upper bounds for their Banach-Mazur distance.
\end{abstract}
\maketitle

\section{Introduction}
Given two isomorphic Banach spaces $X$ and $Y$, their Banach-Mazur distance $\bmd(X,Y)$ is defined as the infimum of distortions $\|T\|\cdot\|T^{-1}\|$ taken over all isomorphisms $T:X\to Y$.

The classical Banach spaces $L_\infty[0,1]$ and $\ell_\infty$ are isomorphic --- this was demonstrated by
Pe{\l}czy\'nski \cite{Pe58} by these breath-taking lines:

\vspace{-3ex}

\begin{equation}\label{dm1}
L_\infty \simeq  \ell_\infty\oplus A\cong  \ell_\infty\oplus\ell_\infty\oplus A\simeq \ell_\infty\oplus L_\infty, 
\end{equation}

\vspace{-3ex}

\begin{equation}\label{dm2}
\ell_\infty \simeq  L_\infty\oplus B\cong  L_\infty\oplus L_\infty\oplus B\simeq L_\infty\oplus \ell_\infty. 
\end{equation}

Here $L_\infty=L_\infty[0,1]$, $A,B$ are Banach spaces, $\simeq$ denotes 
linear isomorphism (i.e. linear homeomorphism); see
section \ref{ub} for details. It seems that little is known about the Banach-Mazur distance of these spaces, see
e.g.\ \cite{DDLS16}*{page 131}. For instance, if  we follow directly the Pe{\l}czy\'nski decomposition method
(and equip direct sums with the max-norm)
then we get  only 
\[ \bmd(L_\infty[0,1],\ell_\infty\oplus L_\infty[0,1])\le 9,\quad  \bmd(\ell_\infty,\ell_\infty\oplus L_\infty[0,1])\le 9,\]
so 
$ \bmd(L_\infty[0,1],\ell_\infty)\le 81$
because the Banach-Mazur distance is multiplicative.

The space  $\ell_\infty$ is isometric to $C(\beta\omega)$, while $L_\infty[0,1]$ is isometric to $C(K)$
for some compact space $K$.
Recall that if we consider two compacta $K$ and $L$ and the corresponding 
Banach spaces $C(K), C(L)$ of real-valued continuous functions with the usual supremum norm,
 then $\bmd\big( C(K), C(L)\big)<2$ implies that $K$ and $L$ are homeomorphic and, consequently, $C(K)$ is isometric to $C(L)$. 
This was independently proved by Amir \cite{Am65} and Cambern \cite{Ca67}. The threshold $2$ is sharp --- see Cohen and Chu \cite{CC95} for a discussion of this phenomenon.

There are few pairs of compacta $K$ and $L$ for which $\bmd\big( C(K), C(L)\big)$
is determined. 
A remarkable exception is provided by the results due to
 Candido and Galego \cite{CG13},  Gergont and Piasecki \cite{GP24} and  Malec and Piasecki \cite{MP25} where
 countable compacta were investigated from this perspective.
 
Our first result, presented  in section \ref{flb}, is quite general:
if $K$ is a zero-dimensional compact space without isolated points and $L$ is any compactification of the set of natural numbers
then  $\bmd\big( C(K), C(L)\big)\ge 3+2\sqrt{2}$.
Our method is  developed in the next section  where we prove that $\bmd(L_\infty[0,1],\ell_\infty)> 7.41$.
The exact shape of the constant is a bit overwhelming, see \ref{ev}. 
However, it should be noted that the appearance of algebraic numbers in such estimates is unavoidable.
For instance, if $K$ is a space consisting of a nontrivial converging sequence and its limit then
 $\bmd\big(C(K\times K), C(K)\big)=2+\sqrt{5}$ --- this is a consequence of the results 
  from  \cite{CG13} and  \cite{MP25}.
 
The following upper bound is established in Section~\ref{ub}:
 \[ \bmd\big(L_\infty[0,1], \ell_\infty\big) \leq (3 + \sqrt{2})^2<19.49.\]
 The gap between the lower and the upper bound obtained here remains quite large.
Although the precise numerical values of these constants are not crucial from the broader perspective of the isomorphic theory of Banach spaces, obtaining meaningful lower bounds for Banach-Mazur distances
seems to require a deeper understanding of the structure of isomorphisms between the spaces involved.

Our starting point can be described as follows: If $T: C(K)\to C(L)$ is a norm-increasing isomorphism (see section \ref{p}) and
we seek the information for the value of $\|T\|$ then  we focus instead on examining the supremum of
the norms of  $\nu_y\in C(K)^\ast$, where $\nu_y=T^\ast\delta_y$ is the measure on $K$
corresponding to the Dirac measure at  $y\in L$.

\section{Preliminaries}\label{p}

Throughout, $K$ and $L$ denote compact Hausdorff spaces.
Every Banach space of the form $C(K)$ is equipped with the supremum norm and
the dual space $C(K)^\ast$ is identified with the space $M(K)$ of all regular signed Borel measures on $K$ of finite variation. 
Given $f\in C(K)$ and $\mu\in M(K)$, we simply write $\mu(f)$ for $\int_K f \diff\mu$.
Given $\mu\in C(K)^\ast$, its norm $\|\mu\|$ is equal to $|\mu|(K)$. Here $|\mu|$ denotes the variation of $\mu$,
so that if $\mu=\mu^+-\mu^-$ is represented as the difference of two orthogonal nonnegative measures
then $|\mu|=\mu^++\mu^-$.

Both  $\ell_\infty$ and $L_\infty[0,1]$ are isometrically isomorphic to spaces of continuous functions. 
Specifically, $\ell_\infty$ is isometric to  $C(\beta\omega)$, the space of continuous functions
on the  Stone--\v{C}ech compactification of the discrete set $\omega$ of natural numbers, see Semadeni \cite{Se71}*{Chapter IV}.
In turn, $L_\infty[0,1]$ can be isometrically  represented as the space $C(K)$, where $K$ is the Stone space of the measure algebra of the Lebesgue measure; here again Semadeni \cite{Se71} and
Fremlin \cite{Fr3} serve as standard references.
By the measure algebra, we mean here
the quotient algebra $Bor[0,1]/\cN$ where $\cN$ is the $\sigma$-ideal of Lebesgue null sets.
It is well-known that $K$ is a nonseparable extremally disconnected compact space without isolated points.

Note that if we identify $L_\infty[0,1]$ with $C(K)$ as above,  
the Riesz representation theorem gives  $L_\infty^\ast[0,1]=M(K)$.
In what follows, we use a more direct description of the dual space 
$L_\infty^\ast[0,1]$, as established by Yosida and Hewitt \cite{YH}*{Theorem~2.3}: 
every continuous functional on  $L_\infty[0,1]$ can be uniquely represented as a 
finitely additive signed measure $\nu$ on $\text{Bor}[0,1]$ that vanishes 
on $\mathcal{N}$. Furthermore, the norm of the functional is equal to the 
total variation of the corresponding measure, $\|\nu\| = |\nu|([0,1])$.

An operator $T \colon X \to Y$ between Banach spaces is said to be \emph{norm-increasing} 
if $\|x\| \le \|Tx\|$ for every $x \in X$. 
Note that every isomorphism $T$ can be made norm-increasing if we multiply it by the constant $\|T^{-1}\|$.

\begin{remark}\label{p:0}
Given two isomorphic Banach spaces $X$ and $Y$, their Banach--Mazur distance 
$\bmd(X,Y)$ is equal to the infimum of $ \|T\|$,  taken over all 
norm-increasing isomorphisms $T$ from $X$ onto $Y$.
\end{remark}

Suppose that $T:C(K)\to C(L)$ is a norm-increasing isomorphism. 
For each $y\in L$ we denote by $\nu_y$ the signed measure on $K$ defined for $g\in C(K)$ by $\nu_y(g)=Tg(y)$; in other words, $\nu_y=T^\ast\delta_y$. 
In this setting, we note the following.

\begin{lemma}\label{p:1}
The measures $\nu_y$ for $y\in L$ form a 1-norming subset of $M(K)$. 
Moreover, for every $h\in C(L)$ there is $\phi\in C(K)$ such that $\nu_y(\phi)=h(y)$ for every $y\in L$.
\end{lemma}

\begin{proof}
If $g\in C(K)$ and $\|g\|=1$, then $\|Tg\|\ge\|g\|=1$, so there is $y\in L$ such that $|\nu_y(g)|=|Tg(y)|\ge 1$.
This means that $\{\nu_y:y\in L\}$ is a 1-norming set.

For any $h\in C(L)$, there is $\phi \in C(K)$ satisfying $T\phi=h$. 
Then $\nu_y(\phi)=T\phi(y)=h(y)$ for every $y\in L$, as required. 
\end{proof}

The lemma given below  is technical, but it will be the key tool  for establishing our lower bounds.
The statement is general because later it will be used in the context of a regular Borel measure on a compact space but also
for a finitely additive measure on $Bor[0,1]$ (the norm $\|\cdot\|$ discussed here is the supremum norm).

\begin{lemma}\label{var}
Let $\nu$ be a finitely additive signed measure on an algebra $\cA$ of subsets of some space $K$. Suppose also that  we are given
 $A\in\cA$  and $\cA$-measurable functions $\phi, \psi,g$  on $K$.
 
 If $r,t, \eps>0$ are some constants such that:
\begin{enumerate}[(i)]
\item $|\nu|(K) \le t$;
\item $\phi(x)\ge r - \eps$ for every $x\in A$;
\item $\nu(\phi)=0$, $\nu(\psi)=1$; 
\item $\big\| \sigma_1\cdot \phi+\sigma_2\cdot \psi+\sigma_3\cdot g\big\|\le r$ for all $\sigma_1,\sigma_2,\sigma_3\in\{-1,0,1\}$. 
\end{enumerate} 
then $|\nu |(K)\ge 2|\nu(A)|+1/r + (1/r)\cdot |\nu(g)|-(3t/r)\eps$.
\end{lemma} 

\begin{proof}
We can  suppose that the value of $\nu(A)$ is positive (consider $-\nu$ otherwise).
Denote 
\[ \phi_0=\phi\cdot \chi_A, \phi_1=\phi-\phi_0,\quad \psi_0=\psi \cdot \chi_A, \psi_1=\psi-\psi_0,\quad g_0=\chi_A\cdot g, g_1=g-g_0.\]

Note that, by  $(iv)$,  we have  $\|\psi_0\|\le \eps$ and $\|g_0\|\le \eps$ (the latter to be used at the very end). Therefore
\begin{equation}\label{2}
 \nu(\psi_1)=\nu(\psi_0+\psi_1)-\nu(\psi_0)\ge 1-t\eps.
\end{equation} 

Similarly, since  $\| \phi_0 -r\cdot \chi_A\|\le\eps$, we obtain  $\big| \nu(\phi_0) -r\nu(A)\big|\le t\eps$ and hence
\[ \nu(\phi_0)\ge r\nu(A) -t\eps.\]
From this we conclude
\begin{equation}\label{3}
 -\nu(\phi_1)=-\nu(\phi)+\nu(\phi_0)\ge r\nu(A)- t\eps.
\end{equation}
 Suppose that $\nu(g)\ge 0$; then we use $(iv)$ for $\psi_1-\phi_1+g$ (and for $\psi_1-\phi_1-g$ in the negative case).
 Using (\ref{2}) and (\ref{3}) we arrive at
\[ |\nu|(K\sm A)\ge (1/r)\cdot \nu(\psi_1-\phi_1+g_1)\ge (1/r)-(t\eps)/r + \nu(A)- (t\eps)/r+(1/r)\cdot \nu(g_1).\]

Since $\nu(g)\ge \nu(g_1)-t\eps$, we finally get  
\[ |\nu |(K)=|\nu|(A)+|\nu|(K\sm A)\ge 2\nu(A)+1/r + (1/r)\cdot \nu(g)-(3t/r)\eps,\] 
so the proof is complete.
\end{proof}

\section{First lower bound}\label{flb}

Let $L$ be any compactification of $\omega$; in other words, $L$ is a compact space that has a countable dense set of isolated points that is identified with $\omega$.

\begin{theorem}\label{flb:1}
If $K$ is a compact { zero-dimensional} space without isolated points, then
\[ \bmd\big( C(K), C(L) \big) \ge 3+2\sqrt{2} \; (>5.82).\]
\end{theorem}

\begin{proof}
Consider a norm-increasing isomorphism $T: C(K)\to C(L)$ and write $\nu_i=T^\ast\delta_i$ for
$i\in\omega\sub L$.

Let $\Phi$ be a family of functions $\phi=\phi(I,\sigma)\in C(K)$, where

\begin{enumerate}[---]
 \item the set $I\sub \omega$ is finite, 
 \item $\sigma:I\to\{-1,1\}$, 
 \item $T\phi(I,\sigma)(i)=\sigma(i)$ for $i\in I$, while
 \item $T\phi(I,\sigma)(i)=0$ whenever $i\in\omega\sm I$.
\end{enumerate}


Note that the last condition implies that $T\phi(I, \sigma)(x) = 0$ for $x \in L\sm \omega$.
Then let $r=\sup\{\|\phi\|:\phi\in\Phi\}$ and fix $\eps>0$. Since $T$ is norm-increasing, $0 < r \le 1$, and
we shall see in a while that $r<1$ must hold.

We choose $\phi=\phi(I,\sigma)\in\Phi$ such that $\|\phi\|>r-\eps$ and an open set $U\sub K$ such that $\phi|U > r-\eps$ (the negative case will be symmetric).
Using the fact that $K$ has no isolated points, we may pick a sequence of pairwise disjoint nonempty clopen sets $A_k\sub U$. 
Then we examine the functions
\[ g_k=\frac{1}{1+\eps}\cdot \phi +\frac{1-r+2\eps}{1+\eps} \cdot \chi_{A_k}.\]
Here, the coefficients  are chosen so that $\|g_k\|>  1$. On the other hand, there is $k$ such that the value $|\nu_i(A_k)|$ is
 arbitrarily small for every $i\in I$. This means that, by a suitable choice of $k$, we have $|\nu_i(g_k)|<1$ whenever $i\in I$. 
For such $k$ there must be $j\in\omega\sm I$ satisfying $|\nu_j(g_k)|\ge 1$. 
Since $\nu_j(\phi)=0$, we conclude that 
\begin{equation}\label{flb1}
|\nu_j(A_k)|\ge \frac{1+\eps}{1-r+2\eps}.
\end{equation}

Consider $\sigma_j : \{j\} \to \{-1, 1\}$ given by $\sigma_j(j) = 1$.
Then conditions of Lemma \ref{var} are satisfied for $\nu=\nu_j$, $\psi = \phi(\{j\}, \sigma_j) \in \Phi$,  $A=A_k$ and $g=0$.
Using \ref{flb1}, it follows that

\begin{equation}\label{4}
t\ge |\nu|(K)\ge 2\cdot   \frac{1+\eps}{1-r+2\eps}+1/r-3(t/r)\eps.
\end{equation} 

As $\eps>0$ can be arbitrarily small,
\begin{equation}\label{5}
t\ge \frac{2}{1-r}+\frac{1}{r}:= \xi(r).
\end{equation} 

Now it remains to find the minimal value of the function $\xi(\cdot)$ defined on right hand side on $(0,1)$: we have 
\[ \xi'(r)=2/(1-r)^2- 1/r^2\]
so $\xi '(r)=0$ for $r=\sqrt{2}-1$. Hence,
\[ t\ge \xi(\sqrt{2}-1)=3+2\sqrt{2},\]
and we are done.
\end{proof}

\section{Second lower bound}
We expand here the method of the previous section to the specific case of a 
norm-increasing isomorphism $T:L_\infty[0,1]\to \ell_\infty$. Again, our task is to find a lower bound for
$t=\|T\|$.

We adapt the beginning of the proof of Theorem \ref{flb:1} to our setting.
Write $\nu_n=T^\ast\delta_n$ for
$n\in\omega\sub L$. Here $\delta_n\in\ell_\infty^\ast$ is given by $\delta_n(x)=x(n)$.
Then $\nu_n\in L_\infty[0,1]^\ast$ so it may be treated as a finitely additive measure on $Bor[0,1]$
that vanishes on Lebesgue-null sets and its variation is bounded by $t$.

Further, let $\phi_n\in L_\infty[0,1]$ be such that $T\phi_n=e_n$, where $e_n$ is the usual `unit' vector in $\ell_\infty$.
Let $\cF$ be a family of functions $\phi=\phi(I,\sigma)\in L_\infty$, where
\[ \phi(I,\sigma)=\sum_{i\in I} \sigma(i)\phi_i,\]
the set $I$ is finite and $\sigma:I\to\{-1,1\}$.
Then, again,  let $r=\sup\{\|\phi\|:\phi\in\cF\}$.

The new ingredient stems from the following.

\begin{lemma}\label{slb:1}
The series $\sum_n | \phi_n|$ converges almost everywhere to a function bounded by $r$ and, consequently,
$\sum_n a_n\phi_n$ converges almost everywhere for every bounded sequence of $a_n$. 
\end{lemma}

\begin{proof}
This is so since for every finite $I\sub\omega$, 
by choosing the signs appropriately, we have 
\[ \sum_{i\in I}|\phi_i(x)|=\phi(I,\sigma)(x)\le \|\phi(I,\sigma)\|\le r\]
for almost all $x\in [0,1]$.
\end{proof}

At this stage, fix $\eps>0$ and  choose $\phi=\phi(I,\sigma)\in\cF$ such that $\|\phi\|>r-\eps$.
Then there is   a Borel set $B$ such that $\phi|B > r-\eps$ (the negative case will be symmetric).

For every Borel nonnegligible set $A\sub B$ we  first take a look at the function 
 \[ g_A=\frac{1}{1+\eps}\cdot \phi +\frac{1-r+2\eps}{1+\eps} \cdot \chi_{A}.\]
 
 Let $a=\sup_n |\nu_n(A)|$; using Lemma \ref{var} we get

\begin{equation}\label{slb2}
t\gtrsim 2a+ \frac{1}{r}. 
\end{equation}

Moreover, the bound $ t\ge \frac{2}{1-r}+\frac{1}{r}$ from the previous section already gives $t>8$ 
whenever $r\le 1/6$ so for the purposes of  our final theorem we can safely suppose that
$r\ge 1/6$. We now fix $\theta$ (see \eqref{slb1} for such a choice):  
\[\theta=\frac{1-r}{2ar+1},\]
consider the function
\[ \Phi_A=\sum_{n\in \omega\sm I}\nu_n(A)\phi_n, \]
and examine
 \[ f_A=\frac{1}{1+\eps}\big( \phi +({1-r+6\eps}) \cdot \chi_{A} - \theta\cdot \Phi_A\big).\]

 \begin{lemma}\label{slb:2}
For every measure–positive set $A\subset B$ we have $\|f_A\|>1$.  
Moreover, there exists such an $A$ for which $|\nu_i(f_A)|<1$ for every $i\in I$.
\end{lemma}

\begin{proof}
First observe that $a\theta \le 3$ whenever $r \geq 1/6$ (as assumed).  
Recall that the function $\sum_{n\in\omega\setminus I} \varphi_n$ is bounded by $\varepsilon$ on $A$.  
Hence $\theta\Phi_A$ is bounded by $3\varepsilon$ on $A$.  
Thus for $x\in A$ we have
\[ f_A(x)\ge \frac{1}{1+\eps}\big(r-\eps+1-r+6\eps -\theta\Phi_A(x)\big)\ge\frac{1+2\eps}{1+\eps},\]
so $\|f_A\|>1$.

Note that to verify the second statement it is enough to check that for every $\eps'>0$ there is $A\sub B$ 
such that
\[ |\nu_i(A)|<\eps' \mbox{ and } \big| \nu_i(\Phi_A)\big|<\eps',\]
for every $i\in I$. This follows from the fact that the measure $\mu=\sum_{i\in I} |\nu_i|$ has a finite variation
 and the observation  that $\Phi_{A_1}+\Phi_{A_2}=\Phi_{A_1\cup A_2}$ for $A_1\cap A_2=\emptyset$.
It follows that if we divide $A$ into large number of pieces $A_k$ then $\mu(\Phi_{A_k})$ can be made arbitrarily small.
\end{proof}

It follows that for $A$ as in Lemma \ref{slb:2}  there is $n\in\omega\sm I$ such that $|\nu_n(f_A)|\ge 1$. 
Since $\nu_n(\phi)=0$, this implies

\[ \Big| (1-r)\nu_n(A)-\theta \nu_n\big(\Phi_A\big)\Big|\gtrsim 1. \]

Let $\Psi=\sum_{k\neq n} \nu_k(A)\phi_k$; we have $\nu_n(\Phi_A)=\nu_n(A)+\nu_n(\Psi)$ so

\[ v=\Big| (1-r-\theta)\nu_n(A)-\theta \nu_n\big(\Psi)\Big|\gtrsim 1. \]

Here it is convenient to stop keeping precise control of  certain quantities.  
We  write $v \gtrsim 1$ to mean that for every $\eta>0$ one can ensure $v > 1-\eta$ by an appropriate 
choice of the parameters involved in the construction (such as $\varepsilon$, the set $A$, the index $n$, etc.).  
This convention will be used throughout the sequel.

Recall that $\theta=\frac{1-r}{2ar+1}$. It is straightforward to check that
\[v = \frac{ar(1-r)}{2ar+1} \left|2 \nu_n(A) - \frac{\nu_n(\Psi)}{ar}\right| \gtrsim 1,\]
therefore, by the triangle inequality, we obtain
\begin{equation}\label{slb1}
2 |\nu_n(A)| + \frac{|\nu_n(\Psi)|}{ar} \gtrsim \frac{2ar+1}{ar(1-r)}.
\end{equation}

\begin{theorem}\label{slb:4}
\[ \bmd(L_\infty[0,1],\ell_\infty)\ge 7.41.\]
\end{theorem}

\begin{proof}
 We continue the analysis begun in this section and use the notation introduced above.

We know that $t\gtrsim 2a+1/r$ - this bound is effective for larger values of $a$. 
A new lower bound is obtained by applying \eqref{slb1} together with Lemma~\ref{var} to $\psi=\varphi_n$ and $g=h=(1/a)\Psi$.
Since $\theta<1$, the function $g$ satisfies assumption~\ref{var}(iv).
Hence
\[ t\gtrsim \frac{2ar+1}{ar(1-r)}+1/r,\]
and this estimate is better for smaller $a$ (e.g.\ $t\ge 8$ if $a=2$ and $r=1/2$). We find the critical value of $a$:
\[  \frac{2ra+1}{ra(1-r)}+1/r=2a+1/r,\]
which gives 
\[ a=\frac{r+\sqrt{(2-r)r}}{2r(1-r)}.\]

Finally $t$ is bounded from below by the minimum of the function
\[ \xi(r)=2\frac{r+\sqrt{(2-r)r}}{2r(1-r)} +1/r, \quad 0<r<1.\]
According to {\sc WolframAlpha} we have $\xi(r)>7.41$ (note that $\xi(1/2)=4+2\sqrt{3}$).
\end{proof}

Let us record the exact value for the function $\xi$ used above:

\begin{equation}\label{ev}
\min_{r\in (0,1)} \xi(r)=\frac{14 + \sqrt[3]{3554 - 66 \sqrt{33}} + \sqrt[3]{ 3554 +
66 \sqrt{33}}} {6}.
\end{equation}

\section{Upper bound}\label{ub}
Let us briefly examine Pe{\l}czy\'nski's argument for
$ L_\infty[0, 1] \simeq\ell_\infty$, outlined in the introduction by \eqref{dm1} and \eqref{dm2}.

Recall that a subspace of a Banach space is $1$-complemented if there exists a norm-one projection onto this subspace.
We say that a Banach space is $1$-injective if it is $1$-complemented in every superspace.

Consider  a $1$-complemented subspace $Y \sub X$; let  $P : X \to Y$ be a norm-one projection onto $Y$
and $A$ be the kernel of $P$.
Then there is an isomorphism between $X$ and $Y \oplus A$ with distortion at most $3$, given by the formula
 $Tx = (3Px, 3/2(x - Px))$ for $x\in X$. Indeed,  we have $\|T\| \leq 3$ and it is not difficult to verify that $\|T^{-1}\| = 1$.
 The constant $3$ is optimal in general: take $Y={\mathbb R}$ as a subspace of the Banach space $c$ of converging sequences
 and use the fact that $\bmd(c,c_0)=3$, which is a particular case of a result due to Gordon \cite{Go70}.
 
Therefore, \eqref{dm1} and \eqref{dm2} give only 
\[ \bmd\big(L_\infty[0,1],L_\infty[0,1]\oplus\ell_\infty\big)\le 9,\quad \bmd\big(\ell_\infty,L_\infty[0,1]\oplus\ell_\infty\big)\le 9,\]
 and the resulting bound is rather poor:  $\bmd\big(\ell_\infty,  L_\infty[0, 1]\big) \leq 81$.

Hence, to get a better estimate we need  to compose the isomorphisms given by the decomposition method 
in a more sophisticated way.  We  prove the following general result.

\begin{theorem}
\label{ns:3 upper}
Assume that $X, Y$ are Banach spaces such that

\begin{enumerate}[(i)]
 \item there are $1$-complemented subspaces $X' \sub X$, $Y' \sub Y$ with $X$ isometric to $Y'$ and $Y$ isometric to $X'$,
 \item $X$ is isometric to $X \oplus X$, and $Y$ is isometric to $Y\oplus Y$.
\end{enumerate}

Then $\bmd(X, Y) \leq (3 + \sqrt{2})^2$.
\end{theorem}

\begin{proof}
Let $P, R$ be norm-one projections, $P : X \to X'$ and $R : Y \to Y'$. 
Recall that $I_X-P$, $I_Y-R$ are projections onto the complements of $X', Y'$; 
denote these complements by $E, F$, respectively, so that $X \simeq X' \oplus E, Y \simeq Y' \oplus F$.

Then we give names to the isometries whose existence is declared by the assumptions:
 \[ \theta : X \to Y', \quad \eta : Y \to X',\]
 \[ \phi = (\phi_1, \phi_2) : X \to X \oplus X,\quad  \psi = (\psi_1, \psi_2) : Y \to Y \oplus Y.\]

Note that $1/(1 + \sqrt{2}) = \sqrt{2} - 1$; these constants become important below.

Using the Pe{\l}czy\'nski decomposition method, we obtain the following three lines of isomorphisms, proving that $X \simeq Y$
\[X \simeq X' \oplus E \cong Y \oplus E \cong Y \oplus Y \oplus E,\]
\[Y \oplus Y \oplus E \cong Y \oplus X' \oplus E \simeq Y \oplus X \simeq Y' \oplus F \oplus X \cong X \oplus F \oplus X,\]
\[X \oplus F \oplus X \cong X \oplus F \cong Y' \oplus F \simeq Y.\]
Let us now write these three isomorphisms by explicit formulas:

\begin{align*}
T &: X \to Y \oplus Y \oplus E, 
&Tx = \Big(\psi_1 \eta^{-1} P x,
\ (\sqrt{2} - 1)\psi_2 \eta^{-1} P x,
\ x - Px\Big),\\
S &: Y \oplus Y \oplus E \to X \oplus X \oplus F, 
&S(y_1, y_2, e) = \Big(\theta^{-1} R y_1,
\ \eta y_2 + e,
\ y_1 - Ry_1\Big),\\
U &: X \oplus X \oplus F \to Y, 
&U(x_1, x_2, f) = \theta\phi^{-1}\Big((1 + \sqrt{2})x_1,
\ x_2\Big) + f,
\end{align*}

For example, to define $T$ we first map $x$ to the pair $(Px, x-Px)$, then use the isometry $\eta^{-1}$ to view $Px \in X'$ as an element of $Y$ and finally we split $\eta^{-1} P x$ using $\psi$.
It is straightforward to calculate the inverses of these operators:

\begin{align*}
T^{-1}&: Y \oplus Y \oplus E \to X, 
&T^{-1}(y_1, y_2, e) = \eta\psi^{-1}\Big(y_1,
\ (1 + \sqrt{2})y_2\Big) + e, \\
S^{-1}&: X \oplus X \oplus F \to Y \oplus Y \oplus E, 
&S^{-1}(x_1, x_2, f) = \Big(\theta x_1 + f,
\ \eta^{-1} P x_2,
\ x_2 - Px_2\Big), \\
U^{-1}&: Y \to X \oplus X \oplus F, 
&U^{-1}y = \Big((\sqrt{2} - 1)\phi_1 \theta^{-1} R y,
\ \phi_2 \theta^{-1} R y,
\ y - Ry\Big).
\end{align*}

Now, we can compose them to obtain $UST$ (we include the intermediate steps to make the calculations easier to verify)
\[STx = \Big(\theta^{-1} R \psi_1 \eta^{-1} P x,
\ (\sqrt{2} - 1) \eta \psi_2 \eta^{-1} P x + x - Px,
\ \psi_1 \eta^{-1} P x - R\psi_1 \eta^{-1} P x\Big)\]
\begin{align*}
USTx = \theta\phi^{-1}&\Big((1 + \sqrt{2})\theta^{-1} R \psi_1 \eta^{-1} P x,
\ (\sqrt{2} - 1) \eta \psi_2 \eta^{-1} P x + x - Px\Big) \\
&+ \psi_1 \eta^{-1} P x - R\psi_1 \eta^{-1} P x.
\end{align*}
and its inverse $T^{-1}S^{-1}U^{-1}$
\[S^{-1}U^{-1} y = \Big((\sqrt{2} - 1)\theta \phi_1 \theta^{-1} R y + y - Ry,
\ \eta^{-1} P \phi_2 \theta^{-1} R y,
\ \phi_2 \theta^{-1} R y - P\phi_2 \theta^{-1} R y \Big),\]
\begin{align*}
T^{-1}S^{-1}U^{-1} y = \eta\psi^{-1}&\Big((\sqrt{2} - 1)\theta \phi_1 \theta^{-1} R y + y - Ry,
\ (1 + \sqrt{2})\eta^{-1} P \phi_2 \theta^{-1} R y\Big) \\
&+ \phi_2 \theta^{-1} R y - P\phi_2 \theta^{-1} R y.
\end{align*}

Recall that all the direct sums of Banach spaces are equipped with the max-norm. 
Recall also that all the operators $ \phi, \psi, \eta, \theta$  have norm $1$, and so do their inverses. It follows that
\[ \|UST\|\le 3 + \sqrt{2}, \quad \|(UST)^{-1}\|\le 3 + \sqrt{2},\]
so the distortion of this isomorphism is bounded by $(3 + \sqrt{2})^2$. 
\end{proof} 

The authors are aware that a slightly better estimate was obtained by Kania and Lewicki \cite{KL26}, by adjusting the constants in the argument above more effectively.

\begin{corollary}\label{final}
$\bmd\big(\ell_\infty, L_\infty[0,1]\big) \leq (3 + \sqrt{2})^2<19.49$.
\end{corollary}

\begin{proof}
We only need to verify the assumptions of Theorem \ref{ns:3 upper}.
It is well known that both $L_\infty[0,1]$ and $\ell_\infty$ are isometric to their squares. 
Moreover, both $L_\infty[0,1]$ and $\ell_\infty$ are $1$-injective; see e.g.\  \cite{DDLS16}*{2.5}
or \cite{AK06}*{section 4.3}.

We can isometrically embed $\ell_\infty$ as a subspace $X$ of $L_\infty[0,1]$: 
take a sequence of pairwise disjoint open intervals $A_n\sub [0,1]$ and send 
$x\in\ell_\infty$ to $\sum_{n\in \omega} x(n)\chi_{A_n}$.
By $1$-injectivity, there is a norm-one projection from $L_\infty[0,1]$ onto $X$ but in fact in this case such a projection
can be effectively defined. 

Finally, it is  a classical fact that $L_\infty[0,1]$ 
has a weak$^*$-separable dual ball and every such a space embeds isometrically into $\ell_\infty$. 
Again by $1$-injectivity, the embedded copy of $L_\infty[0,1]$ must be $1$-complemented in $\ell_\infty$.
\end{proof}

The constant $(3 + \sqrt{2})^2$ given by Corollary \ref{final} does not seem optimal.
It is also worth mentioning that $\ell_\infty$ and $L_\infty[0,1]$ may not  be isomorphic in some models of set theory without the axiom of choice, see V\"ath  \cite{Va98}. 
One can ask whether it is consistent with ZF that, for example, $100 <\bmd(\ell_\infty, L_\infty[0,1]) < \infty$?

We can also use our results to obtain bounds on the Banach--Mazur distances for certain separable Banach spaces.
Recall that for every compact metric space $K$, there is a unique (up to a homeomorphism) compactification $L$ of $\omega$ such that $L \sm \omega$ is homeomorphic to $K$; see \cite{BV23}*{Proposition 4.3.}.
For the Cantor set $2^\omega$, this space $L$ is sometimes (as in \cites{DV24, Ok05}) called the Pe{\l}czy\'nski compactum due to the classical article from 1965 \cite{Pe65}.
Note that $L$ has to be homeomorphic to $L \times \{0, 1\}$.

\begin{corollary}
Let $L$ be the Pe{\l}czy\'nski compactum. Then
\[3 + 2\sqrt{2} \le \bmd\big(C(L), C(2^\omega)\big)\le (3 + \sqrt{2})^2,\]
\end{corollary}

\begin{proof}
The lower bound follows from Theorem \ref{flb:1}, whose assumptions are clearly satisfied.

The assumption that $L$ is homeomorphic to $L \times \{0, 1\}$ implies that $C(L)$ is isometric to $C(L) \oplus C(L)$. 
A similar statement holds for the Cantor set: $C(2^\omega) \cong C(2^\omega) \oplus C(2^\omega)$.

Since $L$ is metrizable, the homeomorphic copy of $2^\omega$ is a retract of $L$. 
It follows that there exists an extension operator of norm $1$ (see \cite{Pe68}*{Section 2}) and thus a $1$-complemented copy of $C(2^\omega)$ in $C(L)$.
It is also a classical fact that every zero-dimensional separable metric space can be embedded into $2^\omega$ (see, e.g. \cite{Ke95}*{Theorem 7.8}). 
Using the same arguments, $C(2^\omega)$ contains a $1$-complemented copy of $C(L)$. 
Now, Theorem \ref{ns:3 upper} yields the upper bound.
\end{proof}

The most interesting case, without a doubt, is the pair of spaces 
$C[0,1]$ and $C(2^\omega)$. The proof that these spaces are isomorphic
is based on a version of the decomposition method, using, in particular, the fact that
$C[0,1]$ is isomorphic to its $c_0$-direct sum (see \cite{AK06}*{Theorem~2.2.3} for details).
Pe{\l}czy\'nski \cite{Pe68} (in his final remarks, see page~73) mentioned the following
estimate:
\begin{equation}\label{pc}
 \bmd\big(C[0,1], C(2^\omega)\big) \le 12.
\end{equation}

We do not know how to verify Pe{\l}czy\'nski's conjecture, nor have we been able to find any
relevant discussion in the literature. To the best of our knowledge, the status of \ref{pc}
therefore remains unclear. 

We would like to thank Frederick Dashiell for interesting correspondence concerning this topic, Benjamin Vejnar for discussions on the Pe{\l}czy\'nski compactum, and James Hagler for his comments on the proof of Lemma \ref{flb:1}.

\bibliography{refs}

\end{document}